\documentclass[12pt,oneside]{amsart}
\usepackage{geometry}                
\usepackage{graphicx}
\usepackage{amssymb}
\usepackage{epstopdf}

\usepackage{amsmath,amsthm,amscd,amssymb}
\usepackage{latexsym}
\usepackage[colorlinks,citecolor=red,pagebackref,hypertexnames=false]{hyperref}
\numberwithin{equation}{section}
\theoremstyle{plain}
\newtheorem{theorem}{Theorem}[section]
\newtheorem{lemma}[theorem]{Lemma}

\theoremstyle{definition}
\newtheorem{definition}[theorem]{Definition}

\theoremstyle{remark}
\newtheorem{remark}[theorem]{Remark}

\newtheorem{case[theorem]}{Case}

\date{October 14, 2018}      
\author{A. Iosevich and K. Taylor}
\address{Department of Mathematics, University of Rochester, Rochester, NY}
\email{iosevich@math.rochester.edu}
\address{Department of Mathematics, The Ohio State University, Columbus, OH}
\email{taylor.2952@osu.edu} 
\title{\parbox{14cm}{\centering{Finite trees inside thin subsets of ${\mathbb R}^ \MakeLowercase{d}$}}}
\begin{document}
\maketitle
\begin{abstract} Bennett, Iosevich and Taylor proved that compact subsets of ${\Bbb R}^d$, $d \ge 2$, of Hausdorff dimensions greater than $\frac{d+1}{2}$ contain chains of arbitrary length with gaps in a non-trivial interval. In this paper we generalize this result to arbitrary tree configurations. \end{abstract}  

\maketitle

\section{Introduction} 

We begin with a seminal result due to Tamar Ziegler, \cite{Z06}, which generalizes an earlier result due to Furstenberg, Katznelson and Weiss \cite{FKW90}. See also \cite{B86}. 

\begin{theorem} \label{z06} Let $E \subset {\Bbb R}^d$, of positive upper Lebesgue density in the sense that 
$$ \limsup_{R \to \infty} \frac{{\mathcal L}^d \{E \cap {[-R,R]}^d \}}{{(2R)}^d}>0, $$ where ${\mathcal L}^d$ denotes the $d$-dimensional Lebesgue measure. Let $E_{\delta}$ denote the $\delta$-neighborhood of $E$. Let $V=\{ {\bf 0}, v^1, v^2, \dots, v^{k-1}\} \subset {\Bbb R}^d$, where $k \ge 2$ is a positive integer. Then there exists $l_0>0$ such that for any $l>l_0$ and any $\delta>0$ there exists $\{x^1, \dots, x^k\} \subset E_{\delta}$ congruent to $lV=\{ {\bf 0}, lv^1, \dots, lv^{k-1}\}$. \end{theorem} 

In particular, this result shows that we can recover every simplex similarity type and sufficiently large scaling inside a subset of ${\Bbb R}^d$ of positive upper Lebesgue density. It is reasonable to wonder whether the assumptions of Theorem \ref{z06} can be weakened, but the following result due to Falconer (\cite{Falc83}) (see also Maga \cite{Mag10}) shows that conclusion may fail even if we replace the upper Lebesgue density condition with the assumption that the set is of dimension $d$. 

\begin{theorem} \label{mag10} (\cite{Mag10}) For any $d \ge 2$ there exists a full dimensional compact set $A \subset {\Bbb R}^d$ such that $A $ does not contain the vertices of any parallelogram. If $d=2$, then given any triple of points $x^1,x^2,x^3$, $x^j \in A$, there exists a full dimensional compact set $A \subset {\Bbb R}^2$ such that $A$ does not contain the vertices of any triangle similar to $\bigtriangleup x^1x^2x^3$. \end{theorem} 

The general question is to study the distance graph with vertices in a compact set of a given Hausdorff dimension. (For more on graph theory, see \cite{BB98}.) More precisely, let $E$ be a compact subset of ${\Bbb R}^d$, $d \ge 2$, and view its points as vertices of a graph where two vertices $x,y$ are connected by an edge if 
$|x-y|=t$, with $| \cdot |$ denoting the Euclidean distance and $t$ a positive real number. Denote the resulting graph by $G_t(E)$. Theorem \ref{mag10} says that if $d=2$ and the Hausdorff dimension of $E$ is equal to $2$, then $G_t(E)$ does not in general contain a triangle. The situation changes in higher dimensions, as was demonstrated by the first listed author of this paper and Bochen Liu in \cite{ILiu16}. 

\begin{theorem} (\cite{ILiu16}) \label{ILiu16} For every $d \ge 4$ there exists $s_0<d$ such that if the Hausdorff dimension of $E$ is $>s_0$, then $E$ contains vertices of an equilateral triangle. \end{theorem} 

\begin{definition} A {\it path} in a graph is a finite or infinite sequence of edges that connect a sequence of distinct vertices. A path of length $k$ connects a sequence of $(k+1)$-vertices, and we refer to this sequence of vertices as a {\it k-chain}.  
 \end{definition} 

Bennett and the two authors of this paper proved in \cite{BIT16} that if the Hausdorff dimension of $E \subset {\Bbb R}^d$, $d \ge 2$, is greater than $\frac{d+1}{2}$, then $G_t(E)$ contains an arbitrarily long path. More generally, they proved the following. 

\begin{theorem} (Theorem 1.7 in \cite{BIT16}) \label{BIT16} Suppose that the Hausdorff dimension of a compact set $E \subset {\Bbb R}^d$, $d \ge 2$, is greater than $\frac{d+1}{2}$. Then for any $k \ge 1$, there exists an open interval $\tilde{I}$ such that for any  ${\{t_i\}}_{i=1}^k \subset \tilde{I}$ there exists a non-degenerate $k$-chain in $E$ with gaps  ${\{t_i\}}_{i=1}^k$. \end{theorem} 

One of the key aspects of the proof of Theorem \ref{BIT16} is the following estimate. 

\begin{theorem} (Theorem 1.8 in \cite{BIT16}) \label{BIT16upper} Suppose that $\mu$ is a compactly supported non-negative Borel measure such that 
\begin{equation} \label{adupper} \mu(B(x,r)) \leq Cr^{s_{\mu}}, \end{equation} where $B(x,r)$ is the ball of radius $r>0$ centered at $x \in {\Bbb R}^d$, for some $s_{\mu}\in(\frac{d+1}{2}, d]$. Then for any $t_1, \dots, t_k>0$ and $\epsilon>0$,
\begin{equation} \label{cbabove} \mu \times \mu \times \dots \times \mu \{(x^1,x^2, \dots, x^{k+1}): t_i -\epsilon \leq |x^{i+1}-x^i| \leq t_i+\epsilon; \ i=1,2, \dots, k \} \leq C\epsilon^k. \end{equation} \end{theorem} 

For the purposes of this paper we are interested in the special case of Theorem \ref{BIT16} where all the $t_i$'s are equal. Our goal is to extend Theorem \ref{BIT16} to more general configurations. 

\begin{definition} A tree is an (undirected) graph in which any two vertices are connected by exactly one path. \end{definition} 

Our main result is the following. 
\begin{theorem} \label{main} Let $E \subset {\Bbb R}^d$, compact of Hausdorff dimension greater than $\frac{d+1}{2}$ and let $T$ be a tree on $k+1$ vertices. Then there exists a non-empty interval $I$ such that for all $t \in I$, $G_t(E)$ contains $T$ as a subgraph. \end{theorem} 

\vskip.125in 

\begin{remark} For an analogous result in sets of positive upper Lebesgue density, see a result due to Lyall and Magyar in \cite{LM18}. \end{remark} 

\vskip.125in 

\section{Proof of Theorem \ref{main}} 

The proof of Theorem \ref{main} is obtained by streamlining and extending the proof of Theorem \ref{BIT16} in a direct and simple way. 

Let $T$ be a graph on $k+1$ vertices. Enumerate the vertices of $T$ and let ${\mathcal E}(T)$ denote the set of pairs $(i,j)$, $i<j$, such that the $i$th vertex is connected with the $j$'th vertex by an edge. Let $\mu$ be a Borel measure supported on $E$ and define 
\begin{equation}\label{family} {\mathcal T}_{T,t}^{\epsilon}(\mu)=\int \dots \int \left\{ \prod_{(i,j) \in {\mathcal E}(T)} \sigma^{\epsilon}_t(x^i-x^j) \right\} d\mu(x^1) \dots d\mu(x^{k+1}).\end{equation}

\vskip.125in  

It is not difficult to see that Theorem \ref{main} would follow from the following estimates: 

\begin{equation} \label{keyupperboundest} {\mathcal T}_{T,t}^{\epsilon}(\mu) \leq C_k, \end{equation} and 

\begin{equation} \label{keylowerboundest} \liminf_{\epsilon \to 0} {\mathcal T}_{T,t}^{\epsilon}(\mu) \ge c_k>0, \end{equation} where $t \in I$, a non-empty interval, $\epsilon>0$ is taken sufficiently small, and both the upper and lower bounds of $c_k$ and $C_k$ respectively hold independently of $t\in I$ and $\epsilon<\epsilon_0$.  We will prove that these estimates hold when the measure $\mu$ is replaced in each variable by the restriction of $\mu$ to an appropriate subset of $E$ of positive $\mu$-measure.  

\vskip.125in 

In the proof of Theorem \ref{BIT16} in \cite{BIT16}, the upper bound was established using the observation that if $Tf=\lambda*f$, where $\lambda$ is a compactly supported measure satisfying $|\widehat{\lambda}(\xi)| \leq C{|\xi|}^{-\alpha}$ for some $\alpha>0$ and $\mu$ is a compactly supported Borel measure satisfying $\mu(B(x,r)) \leq Cr^s$ for some $s>d-\alpha$, then $T$ is a bounded operator from $L^2(\mu)$ to $L^2(\mu)$. The lower bound was established using an inductive procedure generalizing the argument due to the authors of this paper and Mihalis Mourgoglou in \cite{IMT12}. In this paper we streamline the procedure by proving the upper bound and the lower bound at the same time. 

The key feature of our argument is the following calculation. 

\begin{lemma}\label{robustchain} 
Set $$G=G_{t,\epsilon}(1)= \{x\in E: c<\sigma_t^{\epsilon}*\mu(x) <2^{m(1)}\},$$
where $c>0$ and $m(1)\in \mathbb{N}$. 
There exists a non-empty open interval $I$, an $\epsilon_0>0$, and a choice of $c$, $m(1)$, and $\delta>0$ so that  
$$\mu( G_{t,\epsilon} )  >\delta>0$$
whenever $t\in I$ and $0<\epsilon <\epsilon_0$.
\end{lemma}

To prove this result, let $f^t_{\epsilon}=\sigma_t^{\epsilon}* \mu(x).$ It was proved in \cite{BIT16} that there exists a non-empty open interval $I$ and an $\epsilon_0>0$ so that simultaneously the $L^1(\mu)$ norm of $f^t_{\epsilon}$ is uniformly bounded below and the $L^2(\mu)$ norm is bounded above  for all $t\in I$ and $0<\epsilon<\epsilon_0$.  Denote these uniform lower and upper bounds by $0<c$ and $C$ respectively. Let $\epsilon$ and $t$ be such, and set $f=f^t_{\epsilon}.$

Set $G= G_{t,\epsilon}(1) = \{x: c< f(x) < 2^{m(1)}\}$, where $m(1) \in \mathbb{N}$ is to be determined.  Now,
\begin{equation}\label{bam}c< \int f d\mu(x) = \int_{\{f\le c\}} f d\mu(x)   + \int_G f d\mu(x) + \sum_{l=m(1)}^\infty \int_{\{2^l \le f \le 2^{l+1}\}} f d\mu(x).\end{equation}

It is a straight-forward consequence of Chebyshev's inequality and Cauchy-Schwarz that $\mu({\{2^l \le f \le 2^{l+1}\}}) < C2^{-2l}.$ Plugging this into (\ref{bam}) and taking $m(1)$ sufficiently large, it quickly follows that $\mu(G)$ is bounded below away from zero with constants independent of $\epsilon>0$ and $t$.  \\

By induction, using an identical argument to the one above, one can find the following nested sequence of sets of positive $\mu$-measure (where the lower bound on the measure is independent of $t\in I$ and $\epsilon$ small). 
\begin{lemma}\label{repeat} For $j\in \mathbb{N}$, 
 set $$G_{t,\epsilon}(j+1) = \{x\in G_{t,\epsilon}(j): c(j+1) <\sigma_t^{\epsilon}*\mu_j(x) < 2^{m(j)}\},$$
where $\mu_j(x)$ denotes restriction of the measure $\mu$ to the set ${G_{t, \epsilon }(j)}$ and 
$c(j+1) >0$.  
Then there exists numbers $m(j+1) \in \mathbb{N}$, $c(j+1)>0$, and $\delta_{j+1}>0$, so that if $t\in I$ and $0<\epsilon <\epsilon_0$, then 
$$\mu( G_{t,\epsilon}(j+1) )> \delta_{j+1}>0.$$

\end{lemma}

We now demonstrate the pigeon-holing argument that allows us to deduce \eqref{keyupperboundest} and \eqref{keylowerboundest} when $\mu$ in each variable is appropriately restricted.  \\

Fix $k\in \mathbb{N}$, and let $T$ be a tree on $k+1$ vertices.  We say that a vertex is isolated if it is connected to only one other vertex.  
Let $V(1)$ denote the set of isolated vertices of $T$, and let $x^1, \cdots, x^{N(1)}$ denote the collection of vertices who are connected to at least on vertex in $V(1)$. 
 Let $k_1, \cdots, k_{N(1)}$ denote the number of isolated vertices connected to $x^1, \cdots, x^{N(1)}$ respectively.  \\
 
Consider the expression in \eqref{family}.  Integrating in each $v^j \in V(1)$, we replace each of the expressions 
$$\sigma_t^{\epsilon}(x^i - v^j)d\mu(v^j) \,\,\,\,\,\,\text{ by }  \,\,\,\,\,\, \sigma_t^{\epsilon}*\mu(x^i),$$ whenever $x^i$ is connected to $v^j$.    So, if $v^{j_1}, \cdots v^{j_{k_i}} \in V(1)$ are all connected to $x^i$, then we get an expression of the form \begin{equation}\label{mama}\left( \sigma_t^{\epsilon}*\mu(x^i)\right)^{k_i}\end{equation} in the integrand. \\

The next step is to restrict the vertices $x^1, \cdots, x^{N(1)}$ to the set $G(1)$ as in Lemma \ref{robustchain}.  In this way, for each $x^i$, the expression in \eqref{mama} is bounded above and below by positive constants independent of $t\in I$ and $0<\epsilon<\epsilon_0$.  Due to the positivity of the integrand, we can now consider the expression in \eqref{family} with terms of the form in \eqref{mama} removed.  Finally, let $T(2)$ denote the tree with all of the vertices in $V(1)$ removed, and let  ${\mathcal T}_{T,t}^{\epsilon}(2)$ denote the expression in \eqref{mama} with the above mentioned modifications (so any evidence of the vertices in $V(1)$ has been removed) . \\

We repeat this process.  For $j\in \mathbb{N}$, let $T(j+1)$ denote the tree that is obtained after repeating this process $j$-times.  Let $V(j+1)$ denote the set of isolated vertices of $T(j+1)$, and let $y^1, \cdots, y^{N(j)}$ denote the collection of vertices who are connected to at least on vertex in $V(j+1)$. 
 Let $K_1, \cdots, K_{N(j)}$ denote the number of isolated vertices connected to $y^1, \cdots, y^{N(j)}$ respectively.  \\
 
 Consider the expression in ${\mathcal T}_{T,t}^{\epsilon}(j+1)$.  Integrating in each $v^j \in V(j+1)$, we replace each of the expressions 
$$\sigma_t^{\epsilon}(y^i - v^j)d\mu_{j}(v^j) \,\,\,\,\,\,\text{ by }  \,\,\,\,\,\, \sigma_t^{\epsilon}*\mu_{j}(x^i),$$ whenever $y^i$ is connected to $v^j$.    So, if $v^{j_1}, \cdots v^{j_{K_i}} \in V(j+1)$ are all connected to $y^i$, then we get an expression of the form \begin{equation}\label{mamaj}\left( \sigma_t^{\epsilon}*\mu_{j}(x^i)\right)^{k_i}\end{equation} in the integrand. \\

The next step is to restrict the vertices $y^1, \cdots, y^{N(j)}$ to the set $G(j+1)$ as in Lemma \ref{repeat}.  In this way, for each $y^i$, the expression in ${\mathcal T}_{T,t}^{\epsilon}(j+1)$ is bounded above and below by positive constants independent of $t\in I$ and $0<\epsilon<\epsilon_0$.  Due to the positivity of the integrand, we can now consider the expression in ${\mathcal T}_{T,t}^{\epsilon}(j+1)$ with terms of the form in \eqref{mamaj} removed.  Finally, let $T(j+2)$ denote the tree with all of the vertices in $V(j+1)$ removed.\\

This procedure terminates after a finite number of steps, and we are left with an expression of the form
$$\iint \sigma_t^{\epsilon}(z^1 - z^2)d\mu_{J_1}(z^1)d\mu_{J_2}(z^2)=  \int \sigma_t^{\epsilon}*\mu_{J_1}(z^2) d\mu_{J_2}(z^2),  $$  where we assume with out loss of generality that $J_2\geq J_1$ (so that $G(J_2) \subset G(J_1)$).  If $J_2>J_1$, then this expression is bounded above and below by positive constants independent of $0<\epsilon< \epsilon_0$ and independent of $t\in I$.  We obtain the same conclusion when $J_2=J_1$ by simply restricting the variable $z^2$ to the set $G(J_3)$ defined above in Lemma \ref{repeat}.


\vskip.125in 

\begin{thebibliography}{8}


\bibitem{BIT16} M. Bennett, A. Iosevich and K. Taylor, {\it Finite chains inside thin subsets of ${\Bbb R}^d$},  Anal. PDE \textbf{9} (2016), no. 3, 597-614 (http://arxiv.org/pdf/1409.2581.pdf). 

\bibitem{BB98} B. Bollobas, {\it Modern Graph Theory}, Springer, New York, (1998). 







\bibitem{B86} J. Bourgain, {\it A Szemeredi type theorem for sets of positive density}, Israel J. Math. \textbf{54} (1986), no. 3, 307-331. 

\bibitem{Falc83} K.J. Falconer, {\it Some problems in measure combinatorial geometry associated with Paul Erd\H os}, http://www.renyi.hu/conferences/erdos100/slides/falconer.pdf

\bibitem{Fal86} K. J. Falconer, {\it On the Hausdorff dimensions of distance sets} Mathematika \textbf{32} (1986) 206-212.


\bibitem{FKW90} H. Furstenberg, Y. Katznelson, and B. Weiss, {\it Ergodic theory and configurations in sets of positive density} Mathematics of Ramsey theory, 184-198, Algorithms Combin., 5, Springer, Berlin, (1990).

\bibitem{H62} C. Herz, {\it Fourier transforms related to convex sets}, Ann. of Math. (2) \textbf{75} (1962) 81-92. 

\bibitem{ILiu16} A. Iosevich and B. Liu, {\it Equilateral triangles in subsets of ${\Bbb R}^d$ of large Hausdorff dimension}, (https://arxiv.org/pdf/1603.01907.pdf), Israel Math. J. (accepted for publication), (2016). 

\bibitem{IMT12} A. Iosevich, M. Mourgoglou and K. Taylor, {\it On the Mattila-Sj\"{o}lin theorem for distance sets}, Ann. Acad. Sci. Fenn. Math. \textbf{37}, no.2 , (2012). 





\bibitem{LM18} N. Lyall and A. Magyar, {\it Distance graphs and sets of positive upper density in $\mathbb{R}^d$}, (2018), (arXiv:1803.08916). 

\bibitem{Mag10} P. Maga {\it Full dimensional sets without given patterns}, Real Anal. Exchange, \textbf{36}, 79-90, (2010). 






















\bibitem{Z06} T. Ziegler, {\it Nilfactors of ${\Bbb R}^d$ actions and configurations in sets of positive upper density in ${\Bbb R}^m$}, J. Anal. Math. \textbf{99}, pp. 249-266 (2006). 

\end{thebibliography}
\end{document}